\documentclass[12pt]{article}
\usepackage{amsmath, amsthm, amssymb}
\usepackage{graphicx}
\newtheorem{theorem}{Theorem}[section]
\newtheorem{lemma}[theorem]{Lemma}
\theoremstyle{definition}
\newtheorem{definition}[theorem]{Definition}

\theoremstyle{remark}

\numberwithin{equation}{section}

\begin{document}
\title{Simple Braids for Surface Homeomorphisms}
\author{Kamlesh Parwani\thanks{The author was supported in part by NSF Grant DMS0244529.}}
%\address{Department of Mathematics, University of Houston, Houston, TX 77204, USA}
%\email{parwani@math.uh.edu}
%\subjclass{Primary 37E30, 54H20; Secondary 58F20, 57M60}
%\keywords{Rotation vectors, periodic orbits, surfaces}
\date{September 5, 2003.}

\maketitle

\begin{abstract}
Let $S$ be a compact, oriented surface with negative Euler characteristic and let $f:S\rightarrow S$ be a homeomorphism isotopic to the identity. If there exists a periodic orbit with a non-zero rotation vector $(\vec{p},q)$, then there exists a simple braid with the same rotation vector.
\end{abstract}

\section*{Introduction}

There are many famous theorems about periodic points of annulus maps. The
Poincar\'{e}-Birkhoff theorem asserts that any area-preserving homeomorphism
of the annulus isotopic to the identity has periodic
orbits with rational rotation numbers in the rotation interval of the
homeomorphism (see \cite{Brown&Neumann} and \cite{Franks}). With the additional
twist hypothesis, Aubry and Mather proved the existence of periodic
orbits whose radial order is preserved by the map. Such orbits are called
Birkhoff orbits (see \cite{Katok}) or monotone orbits.

This notion of monotone periodic orbits inspired the definition of topologically
monotone orbits in \cite{Boyland}, where Boyland proved that any
homeomorphism of the annulus, isotopic to the identity, that has a periodic
orbit with a non-zero rotation number $\frac{p}{q}$ also has a topologically monotone
periodic orbit with the same rotation number. A topologically monotone
periodic orbit has the property that the isotopy class of the map, keeping the
periodic orbit fixed as a set, is of finite order.

Then the existence of topologically monotone periodic orbits was established
on the torus for smooth maps in \cite{Le Calvez} and for homeomorphisms in 
\cite{Parwani}. It is natural to ask if a similar theorem can be proved on
other surfaces. The main goal of this paper is to describe an analogous
theorem on general orientable surfaces with negative Euler characteristic.

Section 1 introduces rotation vectors and braids for periodic orbits on
orientable surfaces. Section 2 shows that topologically monotone periodic
orbits (trivial braids) are quite ``rare'' on surfaces with negative Euler
characteristic. This motivates the definition of simple braids (which are
not that rare) in Section 3 where the main theorem about the existence of
simple braids, instead of topologically monotone orbits, is proved on
surfaces of negative Euler characteristic.  So simple braids should be considered as an alternative to trivial braids on surfaces with negative Euler characteristic. Also, the main theorem of this paper can be viewed as a Sharkovskii-type forcing result---non-simple braids force the existence of simple braids.

\section{Definitions and Important Results}

\subsection{Rotation vectors on surfaces}

Throughout this paper, $S$ will represent a compact, oriented surface, with or without boundary, and $f:S\rightarrow S$ will be a homeomorphism isotopic to the identity. In this situation, the \textit{rotation vector} for a periodic point is easy to define. Let $x$\ be
a periodic point of least period $n$. There's an isotopy from the identity\
to $f$ and so we get an arc from $x$\ to $f(x)$---denoted $\{x,f(x)\}$---by
following the isotopy. Concatenate the arcs $\{x,f(x)\},\{f(x),f^{2}(x)\}$%
,...,$\{f^{n-1}(x),x\}$\ to obtain a loop $\gamma $. Take the homology
class $[\gamma ]$ of this loop in the surface and divide it by $n$. The 
\textit{rotation vector}\ of the periodic point $x$\ is $\frac{[\gamma ]}{n}$%
. In fact, the rotation vector of every point in the orbit of $x$ is the
same, and so, we can associate the rotation vector $\frac{[\gamma ]}{n}$ to
the entire periodic orbit.

In the example below (Figure 1), the rotation vector is $\left(\frac{1}{7},\frac{-1}{7}
\right) $. The generators chosen are the inner boundary circles endowed with
the orientations given to the boundary. In general, the rotation vector
depends on the choice of generators. We could also depict the vector as a
3-tuple like $(1,-1,7)$. In general, we write $(\vec{p},q)$,
where $\vec{p}$ is the \textit{homology vector} and $q$ is the
period. 
\begin{figure}[ht]
\centerline{\includegraphics[height=2.5in]{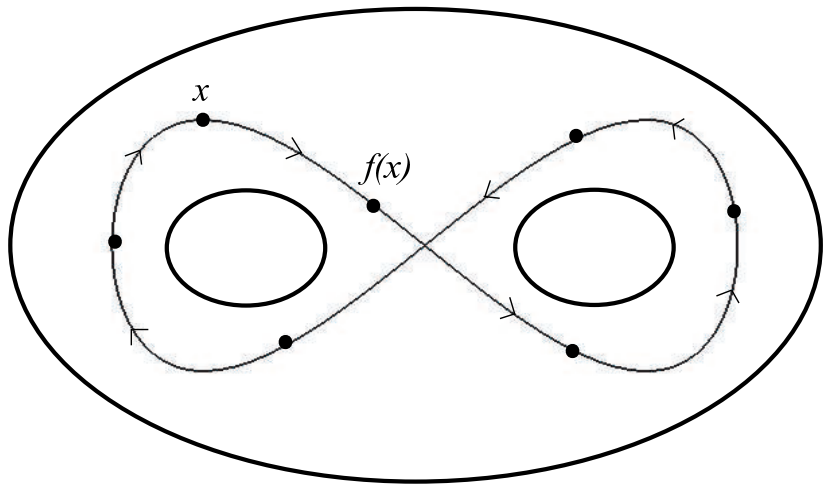}}
\caption{Isotopy loop for a periodic orbit}
\end{figure}

\subsection{The Nielsen-Thurston Classification Theorem and Braids}

Every orientation preserving homeomorphism of an orientable surface with
negative Euler characteristic is isotopic to a homeomorphism $g$ such that
either

a) $g$ is finite order, or

b) $g$ is pseudoAnosov ($pA$), or

c) $g$ is reducible.

A map $g$ is said to be \textit{reducible} if there is a disjoint collection 
$C$ of non-parallel, non-peripheral simple disjoint curves such that $g$
leaves invariant the union of disjoint regular neighborhoods of curves in $C$%
, and the first return map on each complementary component is either of
finite order or $pA$.

This classification theorem was first announced in \cite{Thurston} and the
proofs appeared later in \cite{FLP} and \cite{Casson}.

\smallskip

On any surface, with zero or negative Euler characteristic, we follow Handel
as in \cite{Handel} and examine the isotopy class relative to a periodic
orbit; this will introduce punctures and insure the negative Euler
characteristic required to apply the Nielsen-Thurston Classification
Theorem. When the isotopy class relative to a given periodic orbit is of
finite order, the periodic orbit is called a \textit{finite order periodic
orbit}, and \textit{reducible} and $pA$ \textit{periodic orbits} are defined
similarly. The isotopy class relative to a periodic orbit is also referred
to as the \textit{braid} of the periodic orbit.
\begin{definition}
Let $x$\ and $y$\ be two distinct periodic points of least period $n$\ for
homeomorphisms $f$ and $g$ respectively of the same orientable surface $S$.
Then the orbit of $x$\ ($O(x)$) and the orbit of $y$\ ($O(y)$) have the same 
\textit{braid} if there exists an orientation-preserving homeomorphism $h$\
of $S$ with the property that $h$\ maps $O(x)$\ onto $O(y)$ and the isotopy
class of $h^{-1}fh$\ relative to the orbit of $y$\ is the same as the
isotopy class of $g$\ relative to the orbit of $y$, that is, $%
[h^{-1}fh]_{O(y)}=[f]_{O(y)}$.
\end{definition}

A braid is considered \textbf{trivial} if its isotopy class, relative to the
periodic orbit corresponding to the braid, is of finite order, that is,
there exists a homeomorphism $g$ isotopic to $f$, relative to the periodic
orbit, such that $g^{n}=identity$, for some $n$. And a braid is \textbf{%
non-trivial} if the isotopy class is not of finite order. By definition,
topologically monotone periodic orbits are finite order periodic orbits and
trivial braids, and non-trivial braids are the $pA$ and the reducible
periodic orbits.

Boyland defined a natural partial order $(\vartriangleright )$ on these
braids. If $\alpha $ and $\beta $ are two braids of periodic orbits, then $%
\alpha \vartriangleright \beta $ if and only if the existence of $\alpha $
in any surface homeomorphism $f$ implies the existence of $\beta $ for the
same $f$. The proof of the fact that this is an actual partial order is not
easy and is in \cite{Boyland}.

The existence of topologically monotone periodic orbits on the annulus in 
\cite{Boyland} is established by showing that non-trivial braids of periodic
orbits force the existence of trivial braids of periodic orbits with the
same rotation number (non-trivial $\vartriangleright $ trivial).

\begin{theorem}[Boyland]
Let f be a homeomorphism of the annulus isotopic to the identity. Then if
there exists a periodic orbit with rotation number $\frac{p}{q}$, there
exists a trivial braid with the same rotation number.
\end{theorem}

Then a similar theorem was proved on the torus in \cite{Parwani}.

\begin{theorem}[Parwani]
Let f be a homeomorphism of the torus isotopic to the identity. Then if
there exists a periodic orbit with rotation vector $\left( \frac{p}{q},\frac{%
r}{q}\right) $, there exists a trivial braid with the same rotation vector.
\end{theorem}

The proofs of both theorems rely heavily on the fact that the fundamental
groups of the annulus and the torus are abelian. However, this is not the
case for surfaces with negative Euler characteristic and the situation is
completely different. Section 2 shows that the existence of trivial braids
or finite order periodic orbits cannot always be established. This then
motivates the definition of simple braids in Section 3.

%%%%%%%%%%%%%%%%%%%%%%%%%%%%%%%%%%%%%%%%%%%%%%%%%%%%%%%%%%%%%%%%%%%%%%%

\section{Finite order periodic orbits}

We consider compact, oriented surfaces and homeomorphisms isotopic to the
identity. The main theorem of this section says that, under these
conditions, the isotopy class of any homeomorphism, relative to a periodic
orbit of least period greater than 1, can never be of finite order. This is
why trivial braids cannot exist for periodic orbits with period greater than
1.

The theorem will be proved by a series of simple lemmas. First, some
notation. Let $f:S\rightarrow S$ be a homeomorphism of an oriented surface
that is isotopic to the identity and let $x$ be a periodic point of least
period $n$, $n>1$. Suppose the isotopy class relative to this periodic
orbit, $O(x)$, is of finite order. Let $g$ represent this finite order
class. There exists a hyperbolic metric on this surface for which $g$ is an
isometry (see \cite{Casson}).

\begin{lemma}
For an orientation-preserving isometry g on any connected
oriented surface, if a geodesic segment is fixed (pointwise) by g, then all points are
fixed by g.
\end{lemma}
\begin{proof}
Let $\gamma $ be this geodesic segment that is fixed by $g$.
Assume that $\gamma $ is small, that is, it is contained in a ball with
radius less than the injectivity radius. In this small ball, pick any point $%
x$ and construct a geodesic triangle by picking any two points on $\gamma $
and drawing geodesics from these points to $x$. The two endpoints on $\gamma 
$ are fixed, by assumption, and the orientation preserving isometry sends
geodesics to geodesics and preserves angles and lengths. So the whole
geodesic triangle is fixed, and hence, $x$ is fixed. The point was arbitrary
and so every point is this ball is fixed. Since the surface is connected,
this can be extended to other neighboring balls and eventually to the whole
surface.
\end{proof}

\begin{lemma}
For an orientation-preserving isometry g on any connected
oriented surface, either the fixed points are isolated or g fixes every
point.
\end{lemma}
\begin{proof}
Suppose a fixed point $x$ is not isolated. Then there exists
a fixed point $y$, close enough to $x$, such that there is a unique geodesic
connecting $x$ to $y$. Then the isometry sends this geodesic to itself it fixes every point on
this geodesic. By the lemma above, $g$ fixes every point.
\end{proof}

\begin{lemma}
For an orientation-preserving isometry $g$ on any connected
oriented surface, all the isolated fixed points have index 1.
\end{lemma}
\begin{proof}
The derivative map at the isolated fixed point must preserve
angles and lengths. It is either the identity or some rotation. If it is
the identity, then an entire neighborhood is fixed by $g$ and so every point
on the surface is also fixed. Since we are assuming that the fixed point is
isolated, it must be some rotation. This means that the index is one.
\end{proof}

\begin{theorem}
Let $f:S\rightarrow S$\ be a homeomorphism of a compact, oriented
surface, with negative Euler characteristic, that is isotopic to the
identity and let $x$\ be a periodic point of least period $n$, $n>1$. Then $%
f^{n}$\ is not isotopic to the identity relative to the orbit of $x$.
\end{theorem}
\begin{proof}
Suppose $f^{n}$ is isotopic to the identity relative to the
orbit of $x$. Let $g$ represent this finite order class of $f$. There exists
a hyperbolic metric on this surface for which $g$ is an isometry (see \cite
{Casson}). Clearly $g$ is not the identity as it contains a periodic orbit
of period greater than one, but $g$ is isotopic to the identity because $f$
is. This means that the indices of the fixed points must add up to the Euler
characteristic, which is negative. But by the above lemma, all the fixed
points have index one. So $f$ cannot be isotopic to the identity relative to
the orbit of $x$.
\end{proof}

Here we have established that there are no trivial braids for periodic
orbits with period greater than 1, but what about trivial braids for fixed
points? Does there always exist a finite order fixed point? We will use the following extremely useful description of finite order fixed points which is probably well known. Since the author could not find a reference, a simple proof is provided. First recall that if $f$ is a homeomorphism of surface $S$ isotopic to the identity and if $S$ has negative Euler characteristic, then there exists a unique lift---called the identity lift---of $f$ to the universal cover of $S$ such that the lift commutes with all covering translations. 

\begin{theorem}
Let $f$ be a homeomorphism isotopic to the identity of a surface $S$ with negative Euler characteristic, and let $x$ be a fixed point for $f$. Then $x$ is a finite order fixed point if and only if the identity lift of the map $f$ to the universal cover of $S$ fixes all the lifts of $x$.
\end{theorem}
\begin{proof}
First assume that $S$ has no boundary and consider the following well known exact sequence (see \cite{Birman}).
\begin{equation*}
1 \to \pi_1(S,x) \to MCG(S,x) \to MCG(S) \to 1.
\end{equation*}

Here $MCG(S)$ stands for the mapping class group of $S$ and $MCG(S,x)$ is the mapping class group of $S$ relative to $x$. The homomorphism from $MCG(S,x)$ to $MCG(S)$ is the forgetful homomorphism---one simply removes the restriction that the point $x$ needs to be fixed during the isotopy. The homomorphism from $\pi_1(S,x)$ to $MCG(S,x)$ is called the push homomorphism---given a closed loop based at $x$, we obtain an element in $MCG(S,x)$ by pushing $x$ along the loop.

Now if we restrict our attention to maps which are isotopic to the identity on $S$, then the mapping class group of these homeomorphisms relative to $x$ is isomorphic to $\pi_1(S)$. This isomorphism is easily described. Pick a representative of the element in $MCG(S,x)$ which maps to the identity under the forgetful homomorphism and consider an isotopy of the representative to the identity. Then the loop formed by following $x$ gives an element in $\pi_1(S)$. Equivalently, we may consider the identity lift of the representative and then the covering translation that moves a particular lift of $x$ provides an element in the fundamental group of $S$. It now follows that $f$ is isotopic to the identity relative to $x$ if and only if the identity lift of $f$ fixes all lifts of the point $x$.

Now suppose that $S$ has boundary. If $x$ lies in the interior of $S$, the argument above goes through without any change; here we assume that the mapping class groups are the sets of isotopy classes that fix the boundary of $S$ setwise and not pointwise. If $x$ lies on a boundary circle, $\pi_1(S,x)$ needs to be replaced by the group $\mathbb{Z}$. In other words, $f$ relative to $x$ is isotopic to possibly non-trivial Dehn twists around the boundary circle associated to $x$. The rest of the argument is identical. 
\end{proof}

We answer the question, raised earlier, regarding the existence of finite order fixed points by appealing to a recent result due to Le Calvez (see \cite{Le Calvez2}).

\begin{theorem}[Le Calvez]
Let $F$ be a fixed point free homeomorphism of the plane. Suppose that $F$
commutes with the elements of a discrete group $G$ of orientation-preserving
homeomorphisms that act freely and properly on the plane. Then there exists
a topological foliation of the plane by Brouwer lines that is invariant
under the action of the elements in the group $G$.
\end{theorem}

This allows us to prove the following.

\begin{theorem}
Let $f$ be a homeomorphism, isotopic to the identity, of a compact and oriented surface $S$
with negative Euler characteristic. Then there always exists a finite order
fixed point.
\end{theorem}
\begin{proof}
First suppose that $S$ has no boundary. Let $F$ be the unique lift of $f$ to the universal cover that commutes with all the elements of the fundamental group of $S$. Such a lift does exist since $f$ is isotopic to the identity.
It suffices to show that there exists a fixed point for $F$ because this
would imply that $f$ is isotopic to the identity relative to the fixed point
obtained by projecting the fixed point in the universal cover down to the
surface $S$.

Suppose that $F$ is fixed point free. Then use Le Calvez's result above to obtain a topological foliation of the universal cover, which is homeomorphic to the plane. Since $F$ commutes with all the elements of the fundamental group of $S$, the foliation is invariant under the action of the fundamental group. This means that we obtain a topological foliation of the surface $S$, which is impossible since $S$ has negative Euler characteristic. So $F$ must have a fixed point.

Now suppose $S$ is with boundary. Let $S'$ be another copy of $S$ with the same map $f$ on it and let $\Delta$ (the double) be the surface formed by identifying the corresponding boundary components of $S$ and $S'$. More precisely, every point $x$ in $S$ has a counterpart in $S'$ called $x'$ and $\Delta$ is the surface formed under the identifications $x \sim x'$, for all $x$ on the boundary of $S$. Then a map $g$ on $\Delta$ is defined by $g(x) = f(x)$ for $x \in S$ and $g(x')= (f(x))'$ for $x' \in S'$. Note that $\Delta$ has no boundary and the map $g$ on $\Delta$ is isotopic to the identity. Now the arguments above show that there exists a finite order fixed point for $g$ in $\Delta$---call this point $y$. Without loss of generality, assume that $y$ lies in $S$ and consider the identity lift of $g$ to the universal cover of $\Delta$. This lift fixes all lifts of $y$ and it also leaves invariant all lifts of the subsurface $S$. The restriction of the lift of $g$ to any lift of the subsurface $S$ is the identity lift of $f$ to its universal cover since the covering translations associated to a lift of $S$ are a subset of the covering translations on the universal cover of $\Delta$. Then by Theorem 2.5, it follows that $f$ is isotopic to the identity relative to $y$, and so, $y$ is a finite order fixed point.
\end{proof}

Note that this finite order fixed point must have a zero rotation vector. If
its rotation vector were non-zero, the isotopy loop obtained in the calculation of the rotation vector would have to be essential. Then the lift of this isotopy loop would not be a loop in the
universal cover, implying that $F$ moved the fixed point by some covering
translation. But the finite order fixed point is fixed in the universal
cover, and hence, it must have rotation vector equal to zero.

So the above theorem establishes the existence of a fixed point with
a zero rotation vector for any homeomorphism isotopic to the identity on any
surface with negative Euler characteristic. This is a variation on a
theorem by John Franks who in \cite{Franks2} proved the existence of a fixed
point of positive index and a zero rotation vector under stronger hypotheses.

\begin{theorem}[Franks]
If 0 is in the interior of the convex hull of the recurrent rotation vectors
for an area-preserving diffeomorphism $f$ isotopic to the identity and if
the fixed points are isolated, then $f$ has a fixed point of positive index
and a zero rotation vector.
\end{theorem}

Also note that since all finite order fixed points have to also be
fixed when lifted to the universal cover, no fixed
point with a non-zero rotation vector can be of finite order. So trivial braids
can only exist for fixed points with a zero rotation vector. The next
section concentrates on the case of periodic orbits with non-zero rotation
vectors.

%%%%%%%%%%%%%%%%%%%%%%%%%%%%%%%%%%%%%%%%%%%%%%%%%%%%%%%%%%%%%%%%%%%%%%%%

\section{Simple braids on surfaces}

Let $S$ be a compact, oriented surface with negative Euler characteristic and let $
f:S\rightarrow S$ be a homeomorphism isotopic to the identity. We can always
lift this map to the universal covering space, which in the Poincar\'{e}
disc model can be identified with the interior of the unit disc if $S$ has no boundary. If $S$ has boundary, its boundary components lift to geodesic arcs. In either case, the covering space can be compactified and any lift can be considered as a homeomorphism of closed disc $D$. We will focus on the identity lift $F$ which is the unique lift of $f$ that commutes with all covering translations.
\begin{figure}[ht]
\centerline{\includegraphics[height=2.8in]{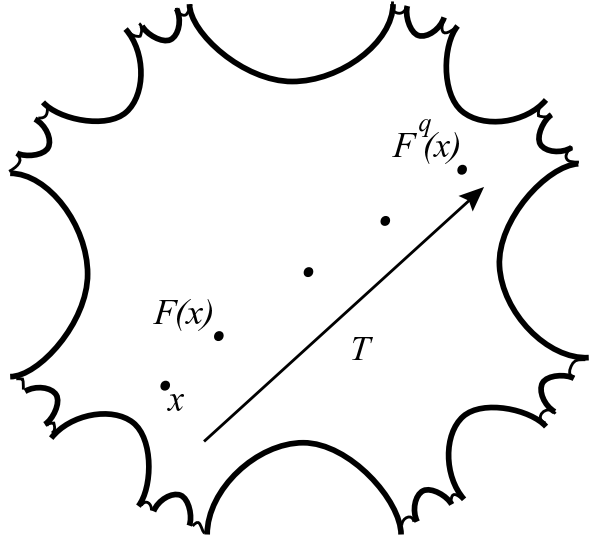}}
\caption{Lift of a surface with boundary}
\end{figure}

First consider the lift of any periodic point in the interior of $S$ with rotation vector $(
\vec{p},q)$ to $D$; assume that $\vec{p}$ is not zero.
If $x$ lies on the lifted orbit, then $F^{q}(x)=T(x)$, for some covering
translation $T$. $D/T$ is homeomorphic to the closed annulus ($A$) and since $F$
commutes with $T$, a homeomorphism $F_{A}$ is induced on $D/T\thickapprox
A$. The lift of the periodic orbit on $D$ projects down to a periodic orbit
on $A$, and since $F^{q}(x)=T(x)$, it has rotation number $\frac{1}{q}$. 

Define a \textit{simple braid} as the periodic orbit that can be lifted, in
the manner described above, to a finite order periodic orbit on the annulus.
So simple braids are trivial when lifted to the annulus. Any periodic orbit on a boundary circle is considered to be a simple braid. In fact, any periodic orbit on a circle trivially has its radial order on the circle preserved under the action of the map, and so we may even call these orbits Birkhoff orbits.

\begin{theorem}
Let $S$\ be a compact, oriented surface with negative Euler characteristic and let $f:S\rightarrow S$\ be a homeomorphism isotopic to the identity. If there exists a periodic orbit with non-zero rotation vector $(\vec{p},q)$, then there exists a simple braid with the same rotation vector.
\end{theorem}
\begin{proof}
If the periodic orbit is already a simple braid, there is nothing to prove. So assume that the orbit lies in the interior of $S$ and then lift the periodic orbit to a periodic orbit with rotation number $\frac{1}{q}$ on the annulus as described above. Boyland's theorem (Theorem
1.2) asserts that there exists a trivial braid with the same rotation vector. If two periodic
points on the annulus have the same rotation vector, then their isotopy
loops, the loops constructed by concatenating isotopy arcs between
successive iterates in order to calculate their rotation numbers, have to be
homotopic (freely). If two loops are homotopic in some covering space, their
projections are also homotopic. So the isotopy loops are homotopic even on $S$. If the isotopy loops are homotopic, they are homologous and the two periodic points have the same rotation vectors.
\end{proof}

\begin{figure}[ht]
\centerline{\includegraphics[height=2.0in]{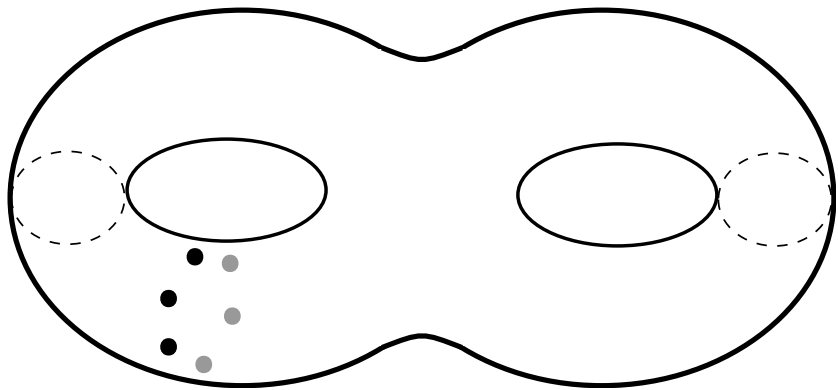}}
\caption{Simple braid on a genus 2 surface}
\end{figure}

Figure 3 depicts a period 6 orbit going around a handle of a genus 2
surface. One is tempted to call this a trivial braid, but by the discussion in Section 2
we know that the periodic orbit cannot be a finite order periodic orbit. The second handle
poses an obstruction to isotoping the map, around the periodic orbit, into a finite order map. Of course, there is no obstruction to doing this on the
torus where trivial braids can exist for all periods (see \cite{Parwani}).

Simple braids are the best we can hope for in general. Even on finite covers
of the surface that one begins with, which have all negative Euler
characteristic, there can be no trivial braids for periodic orbits with
period greater than one. It is only on the compactified, infinite cover, the
annulus, a surface with zero Euler characteristic, where there is some hope
for trivial braids, and that is exactly what happens.

Can this theory of topologically complicated orbits forcing topologically simpler orbits be generalized to the case in which the orbits are not periodic?  This question remains open even on the annulus.

\section*{Acknowledgments}

The author would like to thank John Franks for several useful conversations. This work was supported in part by NSF Grant DMS0244529.

\end{document}